\documentclass[12pt,a4paper]{article}
\usepackage[tbtags]{amsmath}    
\usepackage{amsfonts,amssymb,amsthm,euscript,makeidx,pifont,boxedminipage,multicol,fullpage}
\usepackage{color}

\def\be{\begin{eqnarray}}
\def\ee{\end{eqnarray}}

\def\b*{\begin{eqnarray*}}
\def\e*{\end{eqnarray*}}

\newtheorem{Theorem}{Theorem}[section]
\newtheorem{Lemma}[Theorem]{Lemma}
\newtheorem{Proposition}[Theorem]{Proposition}

\newtheorem{Definition}[Theorem]{Definition}
\newtheorem{Remark}[Theorem]{Remark}

\newcommand{\rmi}{{\rm (i)$\>\>$}}
\newcommand{\rmii}{{\rm (ii)$\>\>$}}
\newcommand{\rmiii}{{\rm (iii)$\>\>$}}
\newcommand{\rmiv}{{\rm (iv)$\>\>$}}

\def \E{\mathbb{E}}
\def \F{\mathbb{F}}

\def \P{\mathbb{P}}
\def \Q{\mathbb{Q}}
\def \R{\mathbb{R}}

\def \N{\mathbb{N}}

\def\Bc{{\cal B}}

\def\Dc{{\cal D}}
\def\Ec{{\cal E}}
\def\Fc{{\cal F}}

\def\Oc{{\cal O}}
\def\Pc{{\cal P}}

\def\Sc{{\cal S}}

\def \Om{\Omega}
\def \om{\omega}
\def \Omb{\overline{\Om}}
\def \omb{\bar{\om}}

\def \eps{\varepsilon}

\def \0{\mathbf{0}}
\def \1{\mathbf{1}}
\def \x{\times}

\def \Fcb{\overline{{\cal F}}}
\def \Fbb{\overline{\F}}

\def \Pcb{\overline{\Pc}}

\def \Ph{\widehat{\P}}

\def \Pb{\overline{\P}}

\def \Bb{\overline{B}}
\def \Gab{\overline{\Gamma}}

\def \SG{\mathrm{SG}}
\def \Ocb{\overline{\Oc}}

\def \Qb{\overline{\Q}}

\def \Ab{\overline{A}}

\title{On the monotonicity principle \\of optimal Skorokhod embedding problem
\footnote{We are grateful to two anonymous referees for their helpful suggestions.
We gratefully acknowledge the financial support of the ERC 321111 Rofirm, the ANR Isotace, the Chairs Financial Risks (Risk Foundation, sponsored by Soci\'et\'e G\'en\'erale), Finance and Sustainable Development (IEF sponsored by EDF and CA).}
}

\author{
	Gaoyue Guo\thanks{CMAP, Ecole Polytechnique, guo@cmap.polytechnique.fr}
	\and Xiaolu Tan\thanks{University of Paris-Dauphine, PSL Research University, CNRS, UMR [7534], CEREMADE.
	tan@ceremade.dauphine.fr}
	\and Nizar Touzi\thanks{CMAP, Ecole Polytechnique, nizar.touzi@polytechnique.edu}
}

\date{\today}

\begin{document}
\bibliographystyle{plain}

\maketitle

\abstract{
	This is a continuation of our accompanying paper \cite{GTT}.
	We provide an alternative proof of the monotonicity principle for the \textit{optimal Skorokhod embedding problem}
	established in Beiglb\"ock, Cox \& Huesmann \cite{BH}.
	Our proof is based on the adaptation of the Monge-Kantorovich duality in our context,
	a delicate application of the optional cross-section theorem,
	and a clever conditioning argument introduced in \cite{BH}.

\vspace{2mm}


\vspace{2mm}


}


\section{Introduction}

	The Skorokhod embedding problem (SEP) consists in constructing a Brownian motion $W$ 
	and a stopping time $\tau$ so that $W_{\tau}$  has some given  distribution.
	Among the numerous solutions of the SEP which appeared in the existing literature, 
	some embeddings enjoy an optimality property with respect to some criterion.
	For instance, the Az\'ema-Yor solution \cite{AY} maximizes the expected running maximum, and 
	Root's embedding \cite{Root} was shown by Rost \cite{Rost} to minimize the expectation of the embedding stopping time.

	Recently, Beiglb\"ock, Cox \& Huesmann \cite{BH} approached this problem by introducing the optimal SEP for some given general criterion.
	Their main result provides a dual formulation in the spirit of optimal transport theory, and a monotonicity principle characterizing optimal embedding stopping times. To the best of our knowledge, 
	all well known solutions to SEP with optimality properties can be interpreted through this unifying principle.
	In addition, the monotonicity principle allows possibly to derive new embeddings with some optimality property as a by-product, see \cite[Sections 2 and 7]{BH} for a more detailed discussion.

	Our main interest in this note is to provide an alternative proof of the last monotonicity principle,
	based on a duality result.
	Our argument follows the classical proof of the monotonicity principle for { the} classical optimal transport problem, see Villani \cite[Chapter 5]{Villani} and the corresponding adaptation by Zaev \cite[Theorem 3.6]{Zaev} for the derivation of the martingale monotonicity principle of Beiglb\"ock \& Juillet \cite{BeiglbockJuillet}. The present continuous-time setting raises however serious technical problems which we overcome in this paper by a crucial use of the optional { cross-section} theorem.

	In the recent literature, there is an important interest in the SEP and the corresponding optimality properties.
	This revival is mainly motivated by its connection to the model-free hedging problem in financial mathematics,
	as initiated by Hobson \cite{Hobson1998}, and { further developed} by many authors \cite{CH, CHO1, CO1, DOR, GTT, Hou, KTZ, OblojSpoida}, etc.

	Finally, we emphasize that the connection between the model-free hedging problem and the optimal transport theory was introduced simultaneously by Beiglb\"ock, Henry-Labord\`ere \& Penkner \cite{BHP} in thediscrete-time { case},
	and Galichon, Henry-Labord\`ere \& Touzi \cite{GHT} in the continuous-time case.
	We also refer to the subsequent literature on martingale optimal transport by \cite{BHT, BBKN14, BT, DS, DS2, HOST, HXT, HT, HK, Hou}, etc.

	In the rest of the paper, we formulate the monotonicity principle in Section \ref{sec:MP}, and then provide our proof in Section \ref{sec:proof}.

\section{Monotonicity principle of optimal Skorokhod embedding problem}
\label{sec:MP}

\subsection{Preliminaries}

	Let $\Om \subset C(\R_+, \R)$ be the canonical space of all continuous functions 
	{ $\om=(\om_t)_{t\ge 0}$} on $\R_+$ such that $\om_0 = 0$,
	{ $B=(B_t)_{t\ge 0}$} denote the canonical process,
	and let $\F = (\Fc_t)_{t \ge 0}$ { be} the canonical filtration generated by $B$.
	Notice that $\Om$ is a Polish space under the compact convergence topology, and its Borel $\sigma-$field
	is given by $\Fc := \bigvee_{t \ge 0} \Fc_t$.
	Denote by $\Pc(\Om)$ the space of all (Borel) probability measures on $\Om$ 
	and by $\P_0\in\Pc(\Om)$ the Wiener measure on $\Om$, under which $B$ is a Brownian motion.

	We next introduce an enlarged canonical space $\Omb := \Om\times\R_+$,
	equipped with canonical element $\Bb:=(B,T)$ defined by
	\b*
		B(\omb)~:=\om ~~\mbox{and}~~ T(\omb)~:=~\theta, ~~~\mbox{for all}~ \omb=(\om,\theta)\in\Omb,
	\e*
	and the canonical filtration $\Fbb=(\Fcb_t)_{t\ge 0}$ defined by 
	\b*
		\Fcb_t &:=& \sigma(B_u, u\le t) \vee \sigma\big(\{T\le u\}, u\le t\big),
	\e*
	so that the canonical variable $T$ is { an} $\Fbb-$stopping time.
	In particular, we have the $\sigma-$field $\Fcb_T$ on $\Omb$.
	Define also $\Fcb^0 := \sigma(B_t, t \ge 0)$ as { the} $\sigma-$field on $\Omb$ generated by $B$.
	Under the product topology, $\Omb$ is still a Polish space, and its Borel $\sigma-$field is given by
	$\Fcb := \bigvee_{t \ge 0} \Fcb_t$.
	Similarly, we denote by $\Pcb(\Omb)$ the set of all (Borel) probability measures on $\Omb$.

	Next, for every $\omb = (\om, \theta) \in\Omb$ and $t\in\R_+$, we define the stopped path by 
	$\om_{t\wedge\cdot}:=\big(\om_{t\wedge u}\big)_{u\ge 0}$
	and $\omb_{t \wedge \cdot} :=(\om_{t \wedge \cdot}, t \wedge \theta)$.
	For every $\omb=(\om,\theta), \omb'=(\om',\theta')\in\Omb$,
	we define the concatenation $\omb\otimes\omb'\in\Omb$ by
	\b*
		\omb\otimes\omb'&:=&(\om\otimes_{\theta}\om', \theta+\theta'),
	\e* 
	where
	\b*
		\big(\om\otimes_{\theta}\om'\big)_t
		&:=&
		\om_t\mathbf{1}_{[0,\theta)}(t)
		~+~
		\big(\om_{\theta}+\om'_{t-\theta}\big)\mathbf{1}_{[\theta,+\infty)}(t),
		~\mbox{for all}~ t\in\R_+.
	\e*
	Let $\xi : \Omb \to \R$ be a non-anticipative $\Fcb-$random variable, i.e. $\xi(\om,\theta)=\xi(\om_{\theta\wedge \cdot},\theta)$ for all $(\om,\theta)\in\Omb$. 
	In the following of the paper, we define, for each $\Pb \in \Pcb(\Omb)$, the expectation $\E[\xi] := \E[ \xi^+ ] - \E[ \xi^-]$, by the convention $\infty - \infty = - \infty$.

\subsection{The optimal Skorokhod embedding problem}

	We now introduce an optimal Skorokhod embedding problem and its dual problem.
	Let $\mu$ be a centered probability measure on $\R$,
	i.e. admitting first order moment and with zero mean,
	we then introduce the set of all embeddings by
	\b*
		\Pcb(\mu)
		&:=&
		\big\{\Pb\in\Pcb: B_T \stackrel{\Pb}{\sim} \mu \big\} ,
	\e*
	with
	\be \label{eq:def_Pcb0}
		\Pcb 
		&:=&
		\big\{ \Pb \in \Pcb(\Omb) ~:
		~B ~\mbox{is an}~ 
		\Fbb-\mbox{Brownian motion and}~ \nonumber \\
		&&~~~~~~~~~~~~~~~~~~~~~~~~~~
		B_{T\wedge\cdot} ~\mbox{is}~\mbox{uniformly integrable under }  \Pb \big\}.
	\ee
	For the given non-anticipative functional $\xi$,
	we define the optimal Skorokhod embedding problem (with respect to $\mu$ and $\xi$) by
	\be \label{eq:P}
		P(\mu)
		&:=&
		\sup_{\Pb\in \Pcb(\mu)}
		\E^{\Pb} 
		\big[ \xi \big].
	\ee

	\begin{Remark}
		The above problem is in fact a weak formulation of the optimal SEP.
		As a strong formulation, one restricts to the class of ``strong'' stopping times, i.e. the stopping times w.r.t. the Brownian filtration.
		Although most of the well known optimal embeddings are ``strong'' stopping times,
		it seems more natural to consider the weak formulation to obtain the general existence of the optimizers, 
		since the set of all weak embeddings is compact under the weak convergence topology.
		Moreover, in some contexts, it is shown that the optimizer is provided by an embedding in ``weak'' sense,
		see e.g. Hobson \& Pedersen \cite{HP}.
		Finally, when $\mu$ has an atom, it seems well known that the two formulations are not equivalent,
		see also Example 2.11 in our accompanying paper \cite{GTT}.
	\end{Remark}

	We next introduce a dual formulation of the above Skorokhod embedding problem \eqref{eq:P}.
	Let $\Lambda$ denote the space of all continuous functions $\lambda : \R \to \R$ of linear growth,
	and define for every $\lambda\in\Lambda$,
	\b*
		\mu(\lambda)&:=&\int_{\R}\lambda(x) \mu(d x).
	\e* 
	Define further
	\b*
		\Dc
		&:=&
		\big\{(\lambda, S)\in\Lambda\x\Sc ~
		: \lambda(\om_{t}) + S_{t}(\om)~\ge~ \xi(\om, t),
		~ \mbox{for all}~ t \ge 0, ~\P_0-\mbox{a.s.}
		\big\},
	\e*
	where $\Sc$ denotes the collection of all  $\F-$strong supermartingales $S=(S_t)_{t\ge 0}$ (which is automatically l\`adl\`ag $\P_0-$a.s.)
	defined on $(\Om, \Fc, \P_0)$ 
	such that $S_0 = 0$ and for some $L>0$,
	\be \label{eq:S_normed}
		\big|S_t (\om)\big|
		&\le &
		L (1 + |\om_t| ),~~ \mbox{for all}~  (\om,t)\in\Omb.
	\ee
	Then the dual problem is given by
	\be \label{eq:D2}
		D(\mu)
		&:=&
		\inf_{(\lambda, S)\in\Dc} ~\mu(\lambda).
	\ee
	
	\begin{Remark}
		By the Doob-Meyer decomposition together with the martingale representation with respect to the Brownian filtration,
		there is some $\F-$predictable process $H=(H_t)_{t\ge 0}$ 
		and non-increasing $\F-$predictable process $A=(A_t)_{t\ge 0}$ ($A_0 = 0$) such that 
		$S_t = (H \cdot B)_t - A_t$ for all $t \ge 0$, $\P_0-$a.s.,
		where $(H \cdot B)$ denotes the stochastic integral of $H$ with respect to $B$ under $\P_0$.
		We then have another dual formulation, by replacing $\Dc$ with
		\b*
			\Dc'
			&:=&
			\big\{(\lambda, H) ~
			: \lambda(\om_{t}) + (H \cdot B)_{t} (\om)~\ge~ \xi(\om, t),
			~ \mbox{for all}~ t \ge 0, ~\P_0-\mbox{a.s.}
			\big\}.
		\e*
		Here, we use the formulation in terms of the set $\Dc$ for ease of presentation.		
	\end{Remark}


\subsection{The monotonicity principle}

	We now introduce the monotonicity principle formulated and proved in Beiglb\"ock, Cox \& Huesmann\cite{BH},
	which provides a geometric characterization of the optimal embedding of problem \eqref{eq:P} in terms of its support.

	Let $\Gab\subseteq\Omb$ be a subset, we define $\Gab^<$ by
	\b*
		\Gab^<
		&:=&
		\big\{\omb=(\om,\theta)\in\Omb
			~:
			\omb_{\theta\wedge\cdot} = \omb'_{\theta \wedge \cdot}
			~\mbox{for some}~
			\omb' \in \Gab
			~\mbox{with}~
			\theta' > \theta
		\big\}.
	\e*
	\begin{Definition} \label{def:SG}
		A pair $(\omb,\omb')\in \Omb\x\Omb$ is said to be a stop-go pair if 
		$\om_{\theta}=\om'_{\theta'}$ and 
		\b*
			\xi(\omb)+\xi(\omb'\otimes \omb'')~>~\xi(\omb\otimes \omb'')+\xi(\omb')
			&\mbox{for all}&
			\omb'' \in\Omb^+,
		\e*
		where $\Omb^{+}:=\big\{\omb=(\om,\theta)\in\Omb: \theta>0\big\}$.
		Denote by $\SG$ the set of all stop-go pairs.
	\end{Definition}

	The following monotonicity principle was introduced and proved in \cite{BH}.

	\begin{Theorem}\label{Th:MP}
		 Let $\xi: \Omb \to \R$ be a Borel non-anticipative random variable.
		Assume that the optimal Skorokhod embedding problem \eqref{eq:P} admits an optimizer 
		 $\Pb^{\ast} \in \Pcb(\mu)$, i.e. $P(\mu)=\E^{\Pb^*}[\xi]$, and the duality $P(\mu) = D(\mu)$ holds.		
		Then there exists a Borel subset $\Gab^{\ast}\subseteq\Omb$ such that
		\b*
			\Pb^* \big[ \Gab^* \big] ~=~ 1
			~&\mbox{and}&~
			\SG ~\cap~\big(\Gab^{\ast<}\x\Gab^{\ast}\big)
			~=~
			\emptyset.
		\e*
	\end{Theorem}

	\begin{Remark}
		\rmi The above monotonicity principle has been formulated and proved in \cite{BH},
		without using the no duality gap condition.

		\noindent \rmii 
		The above duality $P(\mu) = D(\mu)$ has been proved in \cite{BH} (in a slightly stronger formulation) under the condition that $\omb\mapsto \xi(\omb)$ is bounded from above and upper semicontinuous,
		and non-anticipative.

		In our accompanying paper \cite{GTT} (see Theorem 2.4 and Proposition 2.5 of \cite{GTT}),
		we proved the existence of the optimizer $\Pb^*$ as well as the duality $P(\mu) = D(\mu)$
		under the condition that $\xi$ is  non-anticipative, bounded from above and $\theta \mapsto \xi(\om, \theta)$ is upper semicontinuous for $\P_0-$a.e. 
		$\om \in \Om$.
	\end{Remark}

\section{Proof of the main result}
\label{sec:proof}

	Throughout this section, let $\Pb^*$ be an optimizer of problem \eqref{eq:P} in the context of Theorem \ref{Th:MP}.

\subsection{A heuristic proof}

	We start with a purely heuristic argument to illustrate the essential idea in this alternative proof.

	Suppose that there exists a dual minimizer $(\lambda^*,S^*)$ of \eqref{eq:D2}, i.e. 
	\be \label{ineq:dual}
		\lambda^*(\om_t)+S^*_t(\om)~\ge~\xi(\om,t),
		~\mbox{for all}~t \ge 0,
		~\P_0\mbox{-a.s.} ~\mbox{and}~ \mu(\lambda^*)~=~\E^{\Pb^*}[\xi(B,T)],
	\ee
	which implies that $\Gamma:=\{(\om, \theta): \lambda^*(\om_{\theta})+S^*_{\theta}(\om)=\xi(\om,\theta)\}$ has full measure under $\Pb^*$. 
	Assume also for simplicity that $S^*$ is a martingale under $\P_0$.
	We claim that $(\Gamma^<\times\Gamma) \cap \SG=\emptyset$. 
	Otherwise, any pair $(\omb,\omb')\in(\Gamma^<\times\Gamma)\cap \SG$ satisfies the condition
	\b*
		\xi(\omb)+\xi(\omb'\otimes \omb'')~>~\xi(\omb\otimes \omb'')+\xi(\omb')
		&\mbox{for all}&
		\omb'' \in\Omb^+.
	\e*
	Let $\Qb^*_{\omb}$ be the conditional probability of $\Pb^*$ given $\{B_{\theta\wedge \cdot} = \om_{\theta \wedge \cdot}, T > \theta\}$. Then it follows that
	\b*
		\xi(\omb)+\E^{\Qb^*_{\omb}}[\xi(\omb'\otimes \cdot)]&>&\E^{\Qb^*_{\omb}}[\xi(\omb\otimes \cdot)]+\xi(\omb').
	\e*
	On the other hand, notice that the marginal distribution of $\Qb^*_{\omb}$ on $\Om$ is still a Wiener measure.
	Then denoting $(S^*+\lambda^*)(\om, \theta) := S^*_{\theta}(\om) + \lambda^*(\om_{\theta})$,
	one has from \eqref{ineq:dual} that
	\b*
		\xi(\omb)+\E^{\Qb^*_{\omb}}[\xi(\omb'\otimes \cdot)]
		&\le&
		(S^* + \lambda^*)(\omb) + \E^{\Qb^*_{\omb}}[ (S^* + \lambda^*)(\omb' \otimes \cdot)].
	\e*
	Notice that $S^*$ is assumed to be a martingale, and one has from the definition 
	of $\SG$ that $\om_{\theta} = \om'_{\theta'}$, it follows that
	\b*
		(S^* + \lambda^*)(\omb) + \E^{\Qb^*_{\omb}}[ (S^* + \lambda^*)(\omb' \otimes \cdot)]
		=
		\E^{\Qb^*_{\omb}}[(S^* + \lambda^*)(\omb\otimes \cdot)] +  (S^* + \lambda^*)(\omb' ).
	\e*
	Finally, notice that from the definition of $\SG$ and $\Qb^*_{\omb}$, 
	one knows that $\Qb^*_{\omb}[ \omb \otimes \cdot \in \Gamma] = 1$ and $\omb' \in \Gamma$,
	then
	\b*
		\xi(\omb)+\E^{\Qb^*_{\omb}}[\xi(\omb'\otimes \cdot)]
		&\le&
		\E^{\Qb^*_{\omb}}[(S^* + \lambda^*)(\omb\otimes \cdot)] +  (S^* + \lambda^*)(\omb' )\\
		&=&
		\E^{\Qb^*_{\omb}}[\xi(\omb\otimes \cdot)]+\xi(\omb'),
	\e*
	which is a contradiction and we hence obtain that $(\Gamma^{<} \x \Gamma) \cap \SG = \emptyset$.

	\begin{Remark}
		The main technical problem in the above heuristic proof arises from the conditional probability $\Qb^*_{\omb}$ of $\Pb^*$ given $\{B_{\theta \wedge \cdot} = \om_{\theta \wedge \cdot}, T>\theta\}$,
		which should be defined w.r.t. a sub-$\sigma-$field in an almost surely way, creating too many $\Pb^*-$null set to control.
	\end{Remark}

\subsection{An enlarged stop-go set}

	Notice that by Definition \ref{def:SG}, 
	the set $\SG$ is a universally measurable set (co-analytic set more precisely),  but not a Borel set a priori.
	To overcome some measurability difficulty, we will consider as in \cite{BH} another set $\SG^* \subset \Omb \x \Omb$,
	which is Borel.

	Recall that $\Pb^*$ is a fixed optimizer of the problem \eqref{eq:P}, 
	then it admits a family of regular conditional probability distributions (r.c.p.d. see e.g. Stroock \& Varadhan \cite{SV}) $(\Pb^*_{\omb})_{\omb \in \Omb}$ with respect to $\Fcb^0 := \sigma(B_t, t \ge 0)$ on $\Omb$.
	Notice that for $\omb = (\om, \theta)$, the measure $\Pb^*_{\omb}$ is independent of $\theta$, 
	we will denote this family by $(\Pb^*_{\om})_{\om \in \Om}$.
	In particular, one has $\Pb^*_{\om} [B_{\cdot} = \om] = 1$ for all $\om  \in \Om$.
	Next, for every $\omb \in \Omb$, define a probability $\Qb^1_{\omb}$ on $(\Omb, \Fcb)$ by
	\be \label{eq:defQ1}
		\Qb^1_{\omb} [\Ab]
		~:=~
		\int_{\Om} \Pb^*_{\om \otimes_{\theta} \om' }(\Ab) ~\P_0 (d \om'),
		&\mbox{for all}&
		\Ab \in \Fcb.
	\ee
	Intuitively, $\Qb^1_{\omb}$ is the conditional probability with respect to the event 
	$\{ B_{\cdot \wedge \theta} = \om_{\cdot \wedge \theta} \}$.
	We next define, for every $\omb \in \Omb$, a probability $\Qb^2_{\omb}$ by
	\be \label{eq:defQ2}
		\Qb^2_{\omb} [\Ab] 
		~:=~
		\Qb^1_{\omb} \big[ \Ab \big| T > \theta \big] \1_{\{ \Qb^1_{\omb} [ T > \theta] > 0 \}}
		+
		\P_0^{\theta,\om} \otimes \delta_{\{\theta\}} [\Ab]
		\1_{\{ \Qb^1_{\omb} [ T > \theta] = 0 \}},
	\ee
	for all $\Ab \in \Fcb$,
	where $\P_0^{t,\om}$ is the shifted Wiener measure on $(\Om, \Fc)$ defined by
	\b*
		\P_0^{t,\om} [A] ~:=~ \P_0 \big[  \om \otimes_t B \in A \big],
		~~\mbox{for all}~A  \in \Fc.
	\e*
	
	We finally introduce a shifted probability $\Qb^*_{\omb}$ by
	\b*
		\Qb^*_{\omb} [\Ab]
		&:=&
		\Qb^2_{\omb} \big[ \omb \otimes \Bb \in \Ab \big],
		~~\mbox{for all}~
		\Ab \in \Fcb.
	\e*
	and then define a new set $\SG^*$ by
	\be \label{eq:SG_star}
		\SG^* 
		:=
		\big\{ 
			(\omb, \omb') ~: 
			\om_{\theta} = \om'_{\theta'}, 
			\xi(\omb)  + \E^{\Qb^*_{\omb}} [\xi (\omb' \otimes \cdot)] 
			>
			\E^{\Qb^*_{\omb}} [\xi (\omb \otimes \cdot)]  + \xi(\omb') 
		\big\}.
	\ee

	\begin{Lemma} \label{lemm:measurable_issue}
		\rmi The set $\SG^* \subset \Omb \x \Omb$ defined by \eqref{eq:SG_star} is  $\Fcb_T \otimes \Fcb_T-$measurable.

		\vspace{1mm}

		\noindent \rmii Let $\tau \le T$ be a $\Fbb-$stopping time, then the family $(\Ph_{\omb})_{\omb \in \Omb}$ defined by
		$$
			\Ph_{\omb} 
			~:=~
			 \1_{ \{ \tau(\omb) < \theta \}} \Qb^2_{(\om, \tau(\omb))}
			~+~ 
			\1_{\{ \tau(\omb) = \theta \}} \P^{\tau(\omb),\om}_0 \otimes \delta_{\{\theta\}}
		$$
		is a family of regular conditional probability measures of $\Pb^*$ with respect to $\Fcb_{\tau}$,
		i.e. $\omb \mapsto \Ph_{\omb}$ is $\Fcb_{\tau}-$measurable, and for all bounded $\Fcb-$measurable random variable $\zeta$, 
		one has $\E^{\Pb^*}[ \zeta | \Fcb_{\tau}] (\omb) = \Ph_{\omb} [ \zeta]$ for $\Pb^*-$a.e. $\omb \in \Omb$.
	\end{Lemma}
	\proof \rmi 
	Let us denote $[\om]_t := \om_{t \wedge \cdot}$, $[\theta]_t := \theta \1_{\{ \theta \le t \}} + \infty \1_{\{\theta >  t \}}$
	and $[\omb]_t := ([\om]_t, [\theta]_t)$.
	Then by Lemma A.2 of \cite{GTT}, a process $Y : \R_+ \x \Omb \to \R$ is $\Fbb-$optional if and only if it is $\Bc(\R_+) \otimes \Fcb-$measurable,
	and $Y_t(\omb) = Y_t ([\omb]_t)$.
	Further, using Theorem IV-64 of Dellacherie \& Meyer \cite[Page 122]{DM}, 
	it follows that a random variable $X$ is $\Fcb_T-$measurable if and only if it is $\Fcb-$measurable and $X(\omb) = X([\om]_{\theta}, \theta)$ for all $\omb \in \Omb$.
	
	Next, by the definition of $\Qb^1_{\omb}$, $\Qb^2_{\omb}$ and $\Qb^*_{\omb}$, 
	it is easy to see that $\omb \mapsto \big( \Qb^1_{\omb}, \Qb^2_{\omb}, \Qb^*_{\omb} \big)$ are all $\Fcb-$measurable and satisfies $\Qb^*_{\omb} = \Qb^*_{[\om]_{\theta, \theta}}$ for all $\omb \in \Omb$.
	Then it follows that $\omb \mapsto \Qb^*_{\omb}$ is $\Fcb_T-$measurable, 
	and hence by its definition in \eqref{eq:SG_star}, $\SG^*$ is $\Fcb_T \otimes \Fcb_T-$measurable.

	\vspace{1mm}
	
	\noindent \rmii Let $\tau \le T$ be an $\Fbb-$stopping time, we claim that
	\be \label{eq:claim_s.t}
		\mbox{there is some}~ \F- \mbox{stopping time}~ \tau_0 ~\mbox{on}~(\Om, \Fc),~\mbox{s.t.}~
		\tau(\omb) ~=~  \tau_0(\om) \wedge \theta.
	\ee
	Moreover, again by Theorem IV-64 of \cite{DM}, we have
	\be \label{eq:Fcb_tau}
		\Fcb_{\tau}  
		~=~ 
		\sigma \big( B_{\tau \wedge t}, ~t \ge 0 \big) 
		~\vee~ 
		\sigma \big( T \1_{\{ \tau = T\}}, ~\{\tau < T\} \big).
	\ee
	
	Let $(\Ph^0_{\omb})_{\omb \in \Omb}$ be a family of regular conditional probability distribution
	(r.c.p.d. see e.g. Stroock \& Varadhan \cite{SV}) of $\Pb^*$ with respect to $\Fcb_{\tau}$,
	which implies that
	$$
		\Ph^0_{\omb} 
		\big[ 
			B_{\tau \wedge \cdot} = \om_{\tau(\omb) \wedge \cdot}
		\big] = 1
		~\mbox{for all}~\omb \in \Omb;
		~~\mbox{and}~
		\Ph^0_{\omb} [ T = \theta] = 1
		~\mbox{for all}~ \omb \in \{\tau = T\}.
	$$
	It follows that for $\Pb^*-$a.e. $\omb \in \{ \tau = T \}$,
	one has $\Ph^0_{\omb} = \Ph_{\omb} := \P_0^{\tau(\omb), \om} \otimes \delta_{\{\theta \}}$.

	We next focus on the event set $\{\tau < T\}$. Recall that $\Pb^*_{\om}$ is a family of r.c.p.d of $\Pb^*$ with respect to $\sigma (B_t, t \ge 0)$ and $\Qb^1_{\omb}$ are defined by \eqref{eq:defQ1}.
	Then $( \Qb^1_{\om,\tau_0(\om)} )_{\omb \in \Omb}$ is a family of conditional probability measures 
	of $\Pb^*$ with respect to $\sigma \big( B_{\tau_0(\om) \wedge t}, ~t \ge 0 \big)$.
	Further, by the representation of $\Fcb_{\tau}$ in \eqref{eq:Fcb_tau},
	it follows that for $\Pb^*-$a.e. $\omb \in \{ \tau < T\}$, one has 
	$\Ph^0_{\omb} = \Qb^2_{\omb}$.
	
	We now prove the claim \eqref{eq:claim_s.t}.
	For every $\om \in \Om$ and $t \in \R_+$, we denote 
	$\Ab_{\om, t} := \{ \omb' \in \Omb~: \om'_{t \wedge \cdot} = \om_{t \wedge \cdot}, ~\theta' > t \}$.
	Then it is clearly that $\Ab_{\om, t}$ is an atom in $\Fcb_t$, i.e. 
	for any set $C \in \Fcb_t$, one has either $\Ab_{\om, t} \in C$ or $\Ab_{\om, t} \cap C = \emptyset$.
	Let $\omb \in \Omb$ such that $\tau(\omb) < \theta$, and $\theta' > \theta$,
	so that $\omb \in \Ab_{\om, t}$ and $(\om, \theta') \in \Ab_{\om, t}$ for every $t < \theta$.
	Let $t_0 : = \tau(\omb)$, then $\omb \in \Ab_{\om, t_0}$, and $\omb \in \{\tau = t_0\} \in \Fcb_{t_0}$, 
	which implies that $(\om, \theta') \in \Ab_{\om, t_0} \subset \{\tau = t_0\}$ since $\Ab_{\om, t_0}$ is an atom in $\Fcb_{t_0}$.
	It follows that $\tau(\om, \theta') = \tau(\omb)$ for all $\theta' > \theta$ and $\omb \in \Omb$ such that $\tau(\omb) < \theta$.
	Notice that for each $t \in \R_+$, $\{ \omb \in \Omb ~: \tau(\omb) \le t\}$ is $\Fcb_t-$measurable,
	then by Doob's functional representation Theorem, there is some Borel measurable function $f: \Om \x (\R_+ \cup \{ \infty \}) \to \R$ such that
	$\1_{\{ \tau(\omb) \le t\}} = f([\om]_t, [\theta]_t)$.
	It follows that for $\theta_0 \in \R_+$, $\{ \om \in \Om ~: \tau(\om, \theta_0) \le t\}$ is $\Fc_t-$measurable,
	and hence $\om \mapsto \tau(\om, \theta_0)$ is a $\F-$stopping time on $(\Om, \Fc)$.
	Then the random variable $\tau_0 : \Om \to \R_+$ defined by $\tau_0(\om) := \sup_{n \in \N} \tau(\om, n)$ is the required $\F-$stopping time of  claim \eqref{eq:claim_s.t}.

	Finally, we notice that by its definition, 
	one has $\omb \mapsto \Ph_{\omb}$ is $\Fcb-$measurable 
	and satisfies $\Ph_{\omb} = \Ph_{[\omb]_{\theta}}$ for all $\omb \in \Omb$.
	Moreover, we have proved that $\Ph^0_{\omb} = \Ph_{\omb}$ for $\Pb^*-$a.e. $\omb \in \Omb$,
	where $(\Ph^0_{\omb})_{\omb \in \Omb}$ is a family of r.c.p.d. of $\Pb^*$ with respect to $\Fcb_{\tau}$.
	Therefore, $(\Ph_{\omb})_{\omb \in \Omb}$
	is a family of conditional probability measures of $\Pb^*$ with respect to $\Fcb_{\tau}$.
	\qed

	\vspace{2mm}

	To prove Theorem \ref{Th:MP}, we will first prove a closely related result as in \cite{BH}.

	\begin{Theorem}\label{Th:MP_v}
		Let $\xi: \Omb \to \R$ be a Borel non-anticipative random variable.
		Assume that the optimal Skorokhod embedding problem \eqref{eq:P} admits an optimizer 
		 $\Pb^{\ast} \in \Pcb(\mu)$, i.e. $P(\mu)=\E^{\Pb^*}[\xi]$, and the duality $P(\mu) = D(\mu)$ holds.		
		Then there exists a Borel subset $\Gab^{\ast}\subseteq\Omb$ such that
		\b*
			\Pb^* \big[ \Gab^* \big] ~=~ 1
			~&\mbox{and}&~
			\SG^* ~\cap~\big(\Gab^{\ast<}\x\Gab^{\ast}\big)
			~=~
			\emptyset.
		\e*
	\end{Theorem}

\subsection{Technical results}

	We first define a projection operator $\Pi_S : \Omb \x \Omb \to \Omb$ by
	\b*
		\Pi_S \big[ \Ab \big] 
		&:=& 
		\big\{ \omb ~: \mbox{there exists some}~ \omb' \in \Omb~\mbox{such that}~(\omb, \omb') \in \Ab \big \}.
	\e*

	\begin{Proposition} \label{prop:Proj_Tau}
		Let the conditions in Theorem \ref{Th:MP} hold true and 
		$\Pb^*$ be the fixed optimizer of the optimal SEP \eqref{eq:P}.
		Then there is some Borel set $\Gab_0^* \subset \Omb$ such that $\Pb^* [\Gab_0^*] = 1$ and
		for all $\Fbb-$stopping time $\tau \le T$,
		one has
		\be \label{eq:Proj_tau}
			\Pb^* \big[
				\tau < T, 
				~\overline B_{\tau \wedge \cdot} 
				\in \Pi_S \big( \SG^* \cap \big( \Omb \x \Gab_0^* \big)  \big)
			\big] 
			~=~ 0.
		\ee
	\end{Proposition}
	\proof
	\rmi Let us start with the duality result $P(\mu) = D(\mu)$ and the dual problem \eqref{eq:D2}. 
	By definition, we may find a minimizing sequence $\{(\lambda^n,S^n)\}_{n\ge 1} \subset \Dc$ , 
	so that $\mu(\lambda_n)\longrightarrow D(\mu) = P(\mu)$ as $n\longrightarrow\infty$.
	Then, there is some $\Gamma_0 \subset \Om$ s.t. $\P_0(\Gamma_0) = 1$ and
	\be \label{eq:def_eta}
		\eta^n(\omb) := \lambda^n(\om_t)+S^n_t(\om)-\xi(\omb) ~\ge~ 0,
		~~\mbox{for all}~\omb \in \Gamma_0 \x \R_+.
	\ee
	Notice that $(S^n_t)_{t \ge 0}$ are all strong supermartingales on $(\Om, \Fc, \P_0)$ satisfying \eqref{eq:S_normed}.
	It is then also a strong supermartingale on $(\Omb, \Fcb, \Pb^*)$ with respect to $\Fbb$.
	It follows that
	\be \label{eq:cvg_eta}
		0~\le~\E^{\Pb^{\ast}}\big[\eta^n\big]
		=
		\E^{\Pb^*} \big[ \lambda^n(B_T) + S^n_T - \xi  \big]
		\le
		\mu(\lambda^n)-P(\mu)~\longrightarrow~ 0 
		~\mbox{as}~ 
		n \longrightarrow \infty.~
	\ee
	Therefore, we can find some $\Gab_0\subseteq\Omb$ such that $\Pb^{\ast}(\Gab_0) =1$,
	and after possibly passing to a subsequence, 
	\b*
		\eta^{n}(\omb)
		\;\longrightarrow\; 0
		~\mbox{as}~ n\longrightarrow\infty,~~
		\mbox{for all}~~\omb\in\Gab_0.
	\e*
	
	Moreover, since $S^n$ can be viewed as a $\Fbb-$strong supermartingale on $(\Omb, \Fcb, \Pb^*)$,
	then there is some Borel set $\Gab_1 \subset \Omb$ such that  $\Pb^*[ \Gab_1]  = 1$,
	and for all $\omb \in \Gab_1$, 
	$\P_0[ \om \otimes_{\theta} B \in \Gamma_0] = 1$,
	and $(S^n_{\theta + t}(\om \otimes_{\theta} \cdot) )_{t \ge 0}$ is a $\P_0-$strong supermartingale.
	Set $\Gab^*_0 := \Gab_0 \cap \Gab_1$, and we next show that $\Gab^*_0$ is the required Borel set.

	\vspace{1mm}
	
	\noindent \rmii
	Let us consider a fixed pair
	\b*
		(\omb,\omb')&\in&\SG^*~\cap~\big(\Omb \x \Gab^{\ast}_0 \big),
	\e*
	and define
	\b*
		\delta(\omb'')
		&:=&
		\left(\xi(\omb)+\xi(\omb'\otimes \omb'')\right)-\left(\xi(\omb\otimes \omb'')+\xi(\omb')\right),
		~\mbox{for all}~
		\omb''\in\Omb.
	\e*
	By the definition of $\SG^*$ \eqref{eq:SG_star}, one has $\om_{\theta}=\om'_{\theta'}$.
	Then using  the definition of $\eta^n$ in \eqref{eq:def_eta}, 
	it follows that for all $\omb'' \in \Omb$,
	\b*
		\delta(\omb'')
		&=& \lambda^n(\om_{\theta})+S^n_{\theta}(\om)-\eta^n(\omb)
			~-~ \big( \lambda^n(\om'_{\theta'})+S^n_{\theta'}(\om')-\eta^n(\omb') \big) \\
		&& +~
		\lambda^n\big(\om'_{\theta'}+\om''_{\theta''}\big)
			+S^n_{\theta'+\theta''}(\om'\otimes_{\theta'}\om'')
			-\eta^n(\omb'\otimes\omb'') \\
		&&-~
		\big( \lambda^n\big(\om_{\theta} +\om''_{\theta''}\big)
			+ S^n_{\theta+\theta''}(\om\otimes_{\theta}\om'')
			- \eta^n(\omb\otimes\omb'')\big) \\
		&=&
		S^n_{\theta}(\om)-\eta^n(\omb)
			~+~
			S^n_{\theta'+\theta''}(\om'\otimes_{\theta'}\om'')-\eta^n(\omb'\otimes\omb'') \\
		&&
			-~ \big( S^n_{\theta'}(\om')-\eta^n(\omb') 
				+S^n_{\theta+\theta''}(\om\otimes_{\theta}\om'')-\eta^n(\omb\otimes\omb'')\big) \\
		&\le&
		\big ( \eta^n(\omb\otimes\omb'')+\eta^n(\omb')\big)
		~-~
		\eta^n(\omb' \otimes\omb'') \\
		&&
		+~\big ( S^n_{\theta'+\theta''}(\om'\otimes_{\theta'}\om'')-S^n_{\theta'}(\om') \big )
		~-~
		\big ( S^n_{\theta+\theta''}(\om\otimes_{\theta}\om'')-S^n_{\theta}(\om) \big ).
	\e*

	\vspace{1mm}

	\noindent \rmiii
	Let $\tau \le T$ be an $\Fbb-$stopping time, 
	and let $ \big(\Ph_{\omb} \big)_{\omb \in \Omb} $ be the r.c.p.d. of 
	$\Pb^*$ with respect to $\Fcb_{\tau}$ introduced in  Lemma \ref{lemm:measurable_issue}.
	Recall that $\Ph^*_{\omb} := \Qb^*_{(\om, \tau(\omb))}$ for all $\omb \in \{\tau < T\}$ is the shifted probability measures.

	By \eqref{eq:cvg_eta}, there is some set $\Gab^1_{\tau}$ such that
	$\Pb^*[ \Gab^1_{\tau} ] = 1$ and
	\be \label{eq:cvg_eta_tau}
		\E^{\Ph_{\omb}} \big[ \eta^n \big]
		~\longrightarrow~ 0,
		~ \mbox{as}~
		n \longrightarrow \infty,
		~\mbox{for all}~
		\omb \in \Gab^1_{\tau}.
	\ee
	Further, \eqref{eq:cvg_eta} implies  that $ 0 \ge \E^{\Pb^*} [ S^n_T ] \to 0$ as $n \to \infty$.
	Then it follows from the strong supermartingale property of $S^n$ that
	\b*
		~S^n_{\tau} - \E^{\Pb^*}  \big[ S^n_T \big| \Fcb_{\tau} \big] \ge 0,~\Pb^*-\mbox{a.s.}
		~\mbox{and}~
		\E^{\Pb^*}  \big[ S^n_{\tau} - \E^{\Pb^*}  \big[ S^n_T \big| \Fcb_{\tau} \big] \big] 
		\le - \E^{\Pb^*} [ S^n_T ] 
		~\to~ 0,
	\e*
	Hence there is some set $\Gab^2_{\tau} \subset \Omb$ such that
	$\Pb^*[ \Gab^2_{\tau} ] = 1$ and for all $\omb \in  \Gab^2_{\tau}$,
	\be \label{eq:cvg_S}
		0 ~\le~ S^n_{\tau}(\omb) - \E^{\Ph_{\omb}}  [ S^n_T ]
		~\longrightarrow~ 0,
		~\mbox{as}~ n \to \infty,~
	\ee
	after possibly taking some subsequence.
	Moreover, by the definition of $\Pcb$ in \eqref{eq:def_Pcb0}, $B$ is a $\Fbb-$Brownian motion and $B_{T \wedge t}$ is uniformly integrable under $\Pb^*$, and the property holds still under the conditional probability measures.
	Then there is some measurable set $\Gab^3_{\tau} \subset \Omb$ such that 
	$\Pb^*[ \Gab^3_{\tau}] = 1$ 
	and for every $\omb \in \Gab^3_{\tau} \cap \{ \tau < T\}$, one has 
	\be \label{eq:rcpd_support}
		\Ph^*_{\omb} \big[ T > 0 \big] > 0,
		~~
		\Ph^*_{\omb} \big[ \omb_{\tau(\omb) \wedge \cdot} \otimes \overline B \in \Gab^*_0 \big] 
		~=~
		1
		&\mbox{and}&
		\Ph^*_{\omb} \in \Pcb.
	\ee
	Set $\Gab^0_{\tau} := \Gab^1_{\tau} \cap \Gab^2_{\tau} \cap \Gab^3_{\tau}$,
	in the rest of this proof, we show that
	\be \label{eq:contradiction}
		\big( ( \Gab^0_{\tau} \cap \{\tau < T \} ) \x \Omb \big) \cap \SG^* \cap \big( \Omb \x \Gab_0^* \big) 
		&=&
		\emptyset,
	\ee
	which justifies \eqref{eq:Proj_tau}.

	\noindent \rmiv We finally prove \eqref{eq:contradiction} by contradiction.
	Let $(\omb, \omb') \in ( \Gab^0_{\tau} \x \Omb ) \cap \SG^* \cap (\Omb \x \Gab_0^*) $.
	Notice that $\omb' = (\om', \theta') \in \Gab^*_0 \subset \Gab_1$ and for some constant $L_n$,
	\b*
		|S^n_{\theta' +T}(\om' \otimes_{\theta'}B)|&\le& L_n\big|1+\om'_{\theta'}+B_T\big|,
	\e*
	it follows by the supermartingale property, together with the Fatou lemma, that
	\b*
		\E^{\Ph^*_{\omb}}\big[S^n_{\theta'+T}(\om'\otimes_{\theta'}B)\big]
		&\le&
		S^n_{\theta'}(\om').
	\e*
	Moreover, one has $\E^{\Ph^*_{\omb}} [ \eta^n(\omb' \otimes \Bb)] \ge 0$.
	Further, using \eqref{eq:rcpd_support} then \eqref{eq:cvg_eta}, \eqref{eq:cvg_eta_tau}  and \eqref{eq:cvg_S}, we obtain that
	\b*
		0 < \E^{\Ph^*_{\omb}} [\delta]
		\le
		\E^{\Ph^*_{\omb}}\Big[\eta^n(\omb\otimes\Bb)\Big]+\eta^n(\omb')
		-
		\E^{\Ph^*_{\omb}}\Big[S^n_{\tau(\omb) +T}(\om\otimes_{\tau(\omb)}B)\Big]+S^n_{\tau(\omb)}(\om)
		\to 0,
	\e*
	as $n \longrightarrow \infty$,
	which is a contradiction, and we hence conclude the proof.
	\qed

\vspace{2mm}

	Suppose that $\Pi_S \big( \SG^* \cap (\Omb \x \Gab_0^*) \big)$ is Borel measurable on $\Omb$, then by Lemma A.2 of \cite{GTT}, the set 
	$$
		\big \{
			(t, \omb, \omb') \in \R_+ \x \Gab^*_0 \x \Gab^*_0 ~: t < \theta,~\mbox{and}~( (\om, t), \omb') \in \SG^*
		\big \}.	$$
	is an $\Fbb-$optional set.
	Using Proposition \ref{prop:Proj_Tau} together with the classical optional cross-section theorem (see e.g. Theorem IV.86 of Dellacherie \& Meyer \cite{DM}),
	it follows immediately that there is some measurable set $\Gab^*_1 \subset \Omb$ 
	such that $\Pb^*(\Gab^*_1) = 1$ and $\Pi_S \big( \SG^* \cap (\Omb \x \Gab_0^*) \big) \cap \Gab_1^{*<} = \emptyset$.
	However, when the set $\SG^* \cap (\Omb \x \Gab_0^*)$ is a Borel set in $\Omb \x \Omb$,
	the projection set $\Pi_S \big( \SG^* \cap (\Omb \x \Gab_0^*) \big)$ is a priori a $\Bc(\Omb)-$analytic set (Definition III.7 of \cite{DM}) in $\Omb$.
	Therefore, we need to adapt the arguments of the optional cross-section theorem to our context.
	
	\vspace{2mm}

	Denote by $\Ocb$ the optional $\sigma-$field with respect to the filtration $\Fbb$ on $\R_+ \x \Omb$.
	Let $E$ be some auxiliary space, $A \subset \R_+ \x \Omb \x E$, we denote 
	$$
		\Pi_2 (A) ~:=~ \{ \omb ~: \mbox{there is some}~ (t, e) \in \R_+ \x E ~\mbox{such that}~ (t, \omb, e) \in A \},
	$$
	and
	$$
		\Pi_{12}(A) ~:=~ \{ (t,\omb) ~:\mbox{there is some}~ e \in E ~\mbox{such that}~ (t, \omb, e) \in A \}.
	$$

	\begin{Proposition} \label{prop:optional_section}
		Let $\Pb$ be an arbitrary probability measure on $(\Omb, \Fcb)$,
		$(E, \Ec)$ be a Lusin measurable space
		\footnote{A measurable space $(E, \Ec)$ is said to be Lusin if it is isomorphic to a Borel subset of a compact metrizable space (Definition III.16 of \cite{DM}).}.
		Suppose that $A \subset \R_+ \x \Omb \x E$ is a $\Ocb \x \Ec-$measurable set.
		Then for every $\eps > 0$, there is some $\Fbb-$stopping time $\tau$ such that 
		$\Pb [\tau < \infty] ~\ge~ \P[ \Pi_2( A) ] - \eps$
		and $(\tau(\omb), \omb) \in \Pi_{12}(A)$ whenever $\omb \in \Omb$ satisfies $\tau(\omb) < \infty$.
	\end{Proposition}
	\proof
	We follow the lines of Theorem IV.84 of \cite{DM}.
	
	\vspace{1mm}
	
	\noindent \rmi Notice that every Lusin space is isomorphic to a Borel subset of $[0,1]$ (see e.g. Theorem III.20 of \cite{DM}),
	we can then suppose without loss of generality that $(E, \Ec) = ([0,1], \Bc([0,1]))$.
	Then the projection set $\Pi_{12}(A)$ is clearly $\Ocb-$analytic in sense of Definition III.7 of \cite{DM}.	
	
	\vspace{1mm}
	
	\noindent \rmii	
	Using the measurable section theorem (Theorem III.44 of \cite{DM}),
	there is $\Fcb-$random variable $R: \Omb \to \R \cup\{\infty \}$ such that 
	$\Pb[R < \infty] = \Pb[\Pi_2(A)]$ and $R(\omb) < \infty \Rightarrow (R(\omb), \omb) \in \Pi_{12}(A)$.
	The variable $R$ is in fact a stopping time with respect to the completed filtration $\Fbb^{\Pb}$ (see e.g. Proposition 2.13 of \cite{EKT}), but not a $\Fbb-$stopping time a priori.
	We then need to modify $R$ following the measure $\nu$ defined on $\Bc(\R_+) \otimes \Fcb$ by
	$$
		\nu(G) ~:=~ \int \1_G(R(\omb), \omb) \1_{\{ R <\infty \}}(\omb) \Pb(d \omb),
		~~~\forall G \in \Bc(\R_+) \otimes \Fcb.
	$$
	
	\noindent \rmiii We continue by following the lines of item (b) in the proof of Theorem IV.84 of \cite{DM}.
	Denote by $\zeta_0$ the set of all intervals $[[\sigma, \tau [[$, with $\sigma \le \tau$ and $\sigma$, $\tau$ are both $\Fbb-$stoping times.
	Denote also by $\zeta$ the closure of $\zeta_0$ under finite union operation,
	then $\zeta$ is a Boolean algebra which generates the optional $\sigma-$field $\Ocb$.
	Moreover, the debut of a set $C \in \zeta_{\delta}$ (the smallest collection containing $\zeta$ and stable under countable intersection) is a.s. equal to an $\Fbb-$stopping time.
	Further, the projection set $\Pi_{12}(A)$ is $\Ocb-$analytic and hence $\Ocb-$universally measurable.
	Therefore, there exists a set $C \in \zeta_{\delta}$ contained in $\Pi_{12}(A)$ such that $\nu(C) \ge \nu(\Pi_{12}(A)) - \eps$.
	Let $\tau_0$ be the $\Fbb-$stopping time, which equals to the debut of $C$, $\Pb-$a.s.,
	then define $\tau := \tau_0 \1_{\{(\tau_0(\omb), \omb) \in C \}}$, which is a new $\Fbb-$stopping time
	since $\{\omb ~: ( \tau_0(\omb), \omb) \in C\} \in \Fcb_{\tau_0}$ by Theorem IV.64 of \cite{DM}.
	We then conclude the proof by the fact that $\tau$ is the required stopping time.
	\qed

	\vspace{2mm}
	
\subsection{Proof of Theorems \ref{Th:MP_v} and \ref{Th:MP}}

	\noindent {\bf Proof of Theorem \ref{Th:MP_v}}.
	Let us define
	\be \label{eq:def_A}
		A :=
		\big \{
			(t, \omb, \omb') \in \R_+ \x \Gab^*_0 \x \Gab^*_0 ~: t < \theta,~\mbox{and}~((\om, t), \omb') \in \SG^*
		\big \}.
	\ee
	Since $\SG^*$ is a $\Fcb_T \otimes \Fcb_T-$measurable set in $\Omb \x \Omb$,
	it follows (see Lemma A.2. of \cite{GTT}) that the set $A$ defined by \eqref{eq:def_A} satisfies the conditions in Proposition \ref{prop:optional_section} with $E = \Omb$.
	
	We next prove that $\Pi_2(A)$ is $\Pb^*-$null set.
	Indeed, if $\Pb^*[ \Pi_2(A) ] > 0$, by Proposition \ref{prop:optional_section},
	there is some $\Fbb-$stopping time $\tau$ such that 
	$(\tau(\omb), \omb) \in \Pi_{12}(A)$  for all $\omb \in \{ \tau < \infty \}$.
	Notice that $\omb \in \{\tau < \infty\}$ implies that  $(\tau(\omb), \omb) \in \Pi_{12}(A)$ and hence $\tau(\omb) < T$ by the definition of set $A$.
	Therefore, one has $\{\tau < \infty\} = \{ \tau < T \}$,
	and hence $\Pb^*[\tau < \infty] = \Pb^*[\tau < T] > 0$.
	Notice further that $(\tau(\omb), \omb) \in \Pi_{12}(A)$ implies that $(\om, \tau(\omb)) \in \Pi_S(\SG^*)$.
	We then have
	\b*
		0 ~<~ \Pb^* [\tau < T] 
		&\le& 
		\Pb^* 
		\big[
			\tau < T,~
			\overline B_{\tau \wedge \cdot} \in \Pi_S \big( \SG^* \cap \big( \Omb \x \Gab_0^* \big)  \big)
		\big].
	\e*
	This is a contradiction to Proposition \ref{prop:Proj_Tau}.

	Since $\Pi_2(A)$ is a $\Pb^*-$null set, we may obtain a Borel set $\Gab_1^* :=\Gab_0^* \setminus \Pi_2(A)$ 
	such that $\Pb^*[\Gab_1^*] = 1$ and 
	$\Pi_S(\SG^*) \cap \Gab_1^{*<} = \emptyset$.
	Therefore, $\Gab^* := \Gab^*_0 \cap \Gab^*_1$ is the required Borel subset of $\Omb$ in Theorem \ref{Th:MP_v}.	
	\qed
	
	\vspace{3mm}

	\noindent {\bf Proof of Theorem \ref{Th:MP}}.
	Let us define an $\Fbb$-optional process $Z : \R_+ \x \Omb$ by
	$$
		Z_t(\omb) ~=~ Z_t(\om, \theta) ~:=~ \1_{\big\{ t < \theta, ~ \Qb^*_{\om, t}[ T > 0] = 0 \big\}}.
	$$
	Let $\tau$ be an arbitrary $\Fbb-$stopping time, then $\Pb^*[T > \tau] = 1$ implies that
	$\Qb^*_{\om, \tau(\omb)}[T > 0] = 1$ for $\Pb^*-$a.e. $\omb \in \Omb$.
	It follows that 
	$$
		Z_{\tau} = 0, ~~\Pb^*\mbox{-a.s. for all}~\Fbb-\mbox{stopping time}~\tau.
	$$
	Using optional cross-section theorem, one has a Borel set $\Gab^*_2 \subset \Omb$ such that
	$\Pb^*[\Gab^*_2] = 1$ and
	$$
		Z_t = 0,~~\mbox{for all}~t \ge 0,~~\mbox{for every}~\omb \in \Gab^*_2.
	$$
	It is clear that $\Qb^*_{\omb}[\Omb^+] > 0$ for every $\omb \in \Omb$, then by their definition,
	$$
		\SG \cap \big( \Gab_2^{*<} \x \Gab_2^* \big) 
		~\subseteq~
		\SG^* \cap \big( \Gab_2^{*<} \x \Gab_2^* \big) .
	$$
	Then one can conclude the proof by setting $\Gab^* := \Gab^*_0 \cap \Gab^*_1 \cap \Gab^*_2$,
	where $\Gab^*_0$ and $\Gab^*_1$ are the same as in the proof of Theorem \ref{Th:MP_v}.
	\qed

	\begin{Remark}
		\rmi Proposition \ref{prop:Proj_Tau} can be compared to Proposition 6.6 of \cite{BH},
		while the proofs are different.
		Our proof of Proposition \ref{prop:Proj_Tau} 
		is in the same spirit of the classical proof for the monotonicity principle 
		of optimal transport problem (see e.g. Chapter 5 of Villani \cite{Villani}), 
		or martingale optimal transport problem (see e.g. Zaev \cite[Theorem 3.6]{Zaev}),
		based on the existence of optimal transport plan and the duality result.

		\vspace{1mm}

		\noindent \rmii Proposition \ref{prop:optional_section} should be compared to the so-called filtered Kellerer Lemma (Proposition 6.7 of \cite{BH}),
		where a key argument in their proof is Choquet's capacity theory.
		Our proof of Proposition \ref{prop:optional_section} uses crucially an optional section theorem, 
		which is based on a measurable section theorem,
		and the latter is also proved in \cite{DM} using Choquet's capacity theory (see also the review in \cite{EKT}).
	\end{Remark}

\end{document}